\documentclass[11pt]{amsart}
\usepackage{float}
\usepackage{graphicx}
\usepackage{amssymb}
\usepackage{textcomp}
\theoremstyle{plain}
\newtheorem{thm}{Theorem}[section]
\newtheorem{lemma}[thm]{Lemma}

\newtheorem{conj}[thm]{Conjecture}

\theoremstyle{definition}

\newtheorem{example}[thm]{Example}

\def\dim{\mathop{\hbox {dim}}\nolimits}

\newcommand{\fra}{\mathfrak{a}}

\newcommand{\frg}{\mathfrak{g}}
\newcommand{\frh}{\mathfrak{h}}

\newcommand{\frk}{\mathfrak{k}}

\newcommand{\frp}{\mathfrak{p}}

\newcommand{\frt}{\mathfrak{t}}

\newcommand{\bbC}{\mathbb{C}}

\newcommand{\bbR}{\mathbb{R}}

\newcommand{\bbZ}{\mathbb{Z}}

\newcommand{\caC}{\mathcal{C}}

\begin{document}

\title{On the Helgason-Johnson bound}

\author{Chao-Ping Dong}
\address[Dong]{School of Mathematical Sciences, Soochow University, Suzhou 215006,
P.~R.~China}
\email{chaopindong@163.com}

\abstract{Let $G$ be a simple non-compact linear Lie group. Let $\pi$ be any irreducible unitary representation of $G$ with infinitesimal character $\Lambda$ whose continuous part is $\nu$. The beautiful Helgason-Jonson bound in 1969 says that the norm of $\nu$ is upper bounded by the norm of $\rho(G)$, which stands for the half sum of the positive roots of $G$. The current paper aims to give a framework to sharpen the Helgason-Johnson bound when $\pi$ is infinite-dimensional. We have explicit results for exceptional Lie groups. Ingredients of the proof include Parathasarathy's Dirac operator inequality, Vogan pencil, and the unitarily small convex hull introduced by Salamanca-Riba and Vogan.}

\endabstract

\subjclass[2010]{Primary 22E46}

\keywords{}

\maketitle


\section{Introduction}

Let $G$ be a connected simple non-compact linear Lie group which is in the Harish-Chandra class \cite{HC}. For convenience, we further assume that $G$ has finite center. Let $\theta$ be the Cartan involution of $G$, and assume that $K:=G^{\theta}$ is a maximal compact subgroup of $G$. Let $\frg_0$ (resp., $\frk_0$) be the Lie algebra of $G$ (resp., $K$). Then on the Lie algebra level, we have the corresponding Cartan decomposition
$$
\frg_0=\frk_0 + \frp_0,
$$
where $\frp_0$ is the $-1$ eigenspace of $\theta$ on $\frg_0$.
The complexification of $\frg_0$ is written as
$$
\frg:=\frg_0\otimes_{\bbR} \bbC
$$
and similar notation will be used for other groups. We adopt the Killing form on $\frg$, and denote by $\|\cdot\|$ the corresponding norm.

Let $T_f$ be a maximal torus of $K$, and let $\fra_{f, 0}$ be the centralizer of $\frt_{f, 0}$ in $\frp_{0}$. Let $A_f=\exp(\fra_{f, 0})$. Then $H_f=T_f A_f$ is called the fundamental Cartan subgroup of $G$. We have root systems $\Delta(\frg, \frh_f)$, $\Delta(\frk, \frk_f)$ and the restricted root system $\Delta(\frg, \frt_f)$. Denote the corresponding Weyl groups by $W(\frg, \frh_f)$, $W(\frk, \frt_f)$ and $W(\frg, \frt_f)$.
We fix a positive root system $\Delta^+(\frk, \frt_f)$ once for all. Let $\rho_c$ be the half sum of the roots in $\Delta^+(\frk, \frt_f)$. Choose a positive root system $\Delta^+(\frg, \frt_f)$ containing $\Delta^+(\frk, \frt_f)$, and denote the half sum of roots in $\Delta^+(\frg, \frt_f)$ by $\rho(G)$.

We will denote a $K$-type (that is, an irreducible representation of $K$) as $E_{\eta}$, where $\eta$ is its highest weight. Similar notation will apply to $\frk$-types. We may just refer to a $K$-type (or $\frk$-type) by its highest weight. Vogan realized the Langlands classification of the irreducible admissible representations algebraically, where lambda norm and lowest $K$-type played key roles, see Chapter 6 of \cite{V81}. In this setting, an irreducible admissible $(\frg, K)$ module $\pi$ will have infinitesimal character
\begin{equation}\label{inf-char}
\Lambda=(\lambda_a(\mu), \nu)\in \frh^*=\frt^*+\fra^*.
\end{equation}
Let us mention the relevant notation: $E_{\mu}$ is a lowest $K$-type of $\pi$; the vector $\lambda_a(\mu)$ will be recalled in \eqref{lambda-a-mu}, its norm is the \emph{lambda norm} of $\mu$; $G(\lambda_a(\mu))$ is the isotropy group at $\lambda_a(\mu)$ for the $G$ action; $H=TA$ is a maximally split $\theta$-stable Cartan subgroup of $G(\lambda_a(\mu))$. Note that the quasi-split group $G(\lambda_a(\mu))$ is still in the Harish-Chandra class, with Cartan involution $\theta|_{G(\lambda_a(\mu))}$.

Now let us state the following beautiful result.

\medskip
\noindent\textbf{The Helgason-Johnson bound.} Let $\pi$ be any irreducible unitary $(\frg, K)$ module whose infinitesimal character is given by \eqref{inf-char}. Then $\|\nu\|\leq \|\rho(G)\|$.
\medskip

Helgason and Johnson obtained the above result for $\pi$ spherical (that is, $\pi$ contains the trivial $K$-type) in 1969 \cite{HJ}. It is not hard to deduce the general case using the Langlands classification, see Theorem 5.2 of Chapter IV of \cite{BW}. Note that the Helgason-Johnson bound can be tight. For instance, when $G$ is complex, the $\nu$ parameter of the trivial representation does have norm $\|\rho(G)\|$.

However, since a typical irreducible unitary representation of $G$ is infinite-dimensional, it is natural to ask the following question:
Can we improve the Helgason-Johnson bound for infinite-dimensional irreducible unitary representations?

After giving a general framework in Section \ref{sec-framework}, we will handle this question for exceptional Lie groups. It is conceivable that the results reported here will be helpful for understanding the unitary duals.

To make the things clean, we postpone the Hermitian symmetric case to the last section. Our result in this case is Theorem \ref{thm-real-exceptional-HS}. Unless stated otherwise, we assume that $\frk$ has no center from now on. Then the $\frk$-module $\frp$ is irreducible and thus has a highest weight which will be denoted by $\beta$. Let $\pi$ be any infinite-dimensional representation of $G$. If the $K$-type $E_{\mu}$ occurs in $\pi$, then by Section 3 of Vogan \cite{V80}, the $K$-types $E_{\mu+m\beta}$ will show up in $\pi$ for all $m\in\bbZ_{\geq 0}$. Let us call the following $K$-types the \textbf{Vogan pencil} starting from $\mu$
\begin{equation}\label{pencil-mu}
{\rm Pencil}(\mu):=\{\mu+m\beta\mid m\in\bbZ_{\geq 0}\}.
\end{equation}
Let us also introduce
\begin{equation}\label{spin-norm-pencil-mu}
{\rm MP}(\mu):=\min\{\|\mu+m\beta\|_{\rm spin}^2\mid m\in\bbZ_{\geq 0}\},
\end{equation}
where $\|\cdot\|_{\rm spin}$ stands for the spin norm \cite{D13}. It will be recalled in Section \ref{sec-lambda-spin-norm}.

\begin{thm}\label{thm-complex-exceptional}
Let $G$ be a complex exceptional Lie group. Let $\pi$ be any infinite-dimensional irreducible unitary $(\frg, K)$ module whose infinitesimal character is given by \eqref{inf-char}. Then we have the following bounds for $\|\nu\|$:
\begin{itemize}
\item[$\bullet$] ${\rm Lie}(G)$ is of type $G_2:$ $\|\nu\|\leq \sqrt{26}$, while $\|\rho(G)\|=\sqrt{56}$.
\item[$\bullet$] ${\rm Lie}(G)$ is of type $F_4:$ $\|\nu\|\leq \sqrt{78}$, while $\|\rho(G)\|=\sqrt{156}$.
\item[$\bullet$] ${\rm Lie}(G)$ is of type $E_6:$ $\|\nu\|\leq \sqrt{170}$, while $\|\rho(G)\|=\sqrt{312}$.
\item[$\bullet$] ${\rm Lie}(G)$ is of type $E_7:$ $\|\nu\|\leq \sqrt{464}$, while $\|\rho(G)\|=\sqrt{798}$.
\item[$\bullet$] ${\rm Lie}(G)$ is of type $E_8:$ $\|\nu\|\leq \sqrt{1522}$, while $\|\rho(G)\|=\sqrt{2480}$.
\end{itemize}
\end{thm}

Note that the unitary dual of complex $G_2$ has been determined by Duflo \cite{Du} in 1979. Note also that the above bound can be tight. For instance, when $G$ is complex $F_4$, the $\nu$ parameter of the minimal representation does have norm $\sqrt{78}$.

\begin{thm}\label{thm-real-exceptional}
Let $G$ be a real exceptional Lie group. Let $\pi$ be any infinite-dimensional irreducible unitary $(\frg, K)$ module whose infinitesimal character is given by \eqref{inf-char}. Then we have the following bounds for $\|\nu\|$:
\begin{itemize}
\item[$\bullet$] ${\rm Lie}(G)=G_{2(2)}:$  $\|\nu\|\leq \sqrt{\frac{51}{8}}$, while $\|\rho(G)\|=\sqrt{14}$.
\item[$\bullet$] ${\rm Lie}(G)={\rm FI}=F_{4(4)}:$ $\|\nu\|\leq \sqrt{\frac{39}{2}}$, while $\|\rho(G)\|=\sqrt{39}$.
\item[$\bullet$] ${\rm Lie}(G)={\rm FII}=F_{4(-20)}:$ $\|\nu\|\leq \frac{9}{2}$, while $\|\rho(G)\|=\sqrt{39}$.
\item[$\bullet$] ${\rm Lie}(G)={\rm EI}=E_{6(6)}:$ $\|\nu\|\leq \sqrt{42}$, while $\|\rho(G)\|=\sqrt{78}$.
\item[$\bullet$] ${\rm Lie}(G)={\rm EII}=E_{6(2)}:$ $\|\nu\|\leq \sqrt{42}$, while $\|\rho(G)\|=\sqrt{78}$.
\item[$\bullet$] ${\rm Lie}(G)={\rm EIV}=E_{6(-26)}:$ $\|\nu\|\leq \sqrt{37}$, while $\|\rho(G)\|=\sqrt{78}$.
\item[$\bullet$] ${\rm Lie}(G)={\rm EV}=E_{7(7)}:$ $\|\nu\|\leq \sqrt{\frac{231}{2}}$, while $\|\rho(G)\|=\sqrt{\frac{399}{2}}$.
\item[$\bullet$] ${\rm Lie}(G)={\rm EVI}=E_{7(-5)}:$ $\|\nu\|\leq \sqrt{\frac{227}{2}}$, while $\|\rho(G)\|=\sqrt{\frac{399}{2}}$.
\item[$\bullet$] ${\rm Lie}(G)={\rm EVIII}=E_{8(8)}:$ $\|\nu\|\leq \sqrt{380}$, while $\|\rho(G)\|=\sqrt{620}$.
\item[$\bullet$] ${\rm Lie}(G)={\rm EIX}=E_{8(-24)}:$ $\|\nu\|\leq \sqrt{\frac{723}{2}}$, while $\|\rho(G)\|=\sqrt{620}$.
\end{itemize}
\end{thm}

Note that the unitary dual of $G$ when ${\rm Lie}(G)$ is $G_{2(2)}$ has been determined by Vogan \cite{V94} in 1994. Note also that the above bound can be tight. For instance, when ${\rm Lie}(G)$ is EI, the $\nu$ parameter of the representation in the penultimate row of Table 11 of \cite{DDY} does have norm $\sqrt{42}$. This is a minimal representation.
It is worth noting that our bound, among the above cases,  is attained at the trivial $K$-type if and only if $G/K$ is \emph{not} quaternionic as in Appendix C of Knapp \cite{Kn}, which traces back to Wolf \cite{Wo}. See Section \ref{sec-real-exc} for details.

On the $p$-adic side, Dirac operators have been introduced by Barbasch, Ciubotaru and Trapa for graded affine Hecke algebras \cite{BCT}. In that setting one also has Dirac inequality, and an analogue of our results on the $p$-adic side has been deduced in Corollary 5.4 of \cite{BCT}. See also Theorem 1.2 (1) of  Ciubotaru \cite{Ci}. Moreover, Theorem 1.2 (3) of \cite{Ci} gives the second non-unitarity gap by understanding the minimal representations. It is conceivable that there should be an analogue on the real side. We will address this problem in future.

Throughout this paper, the root systems will be adopted as in Appendix C of Knapp \cite{Kn}.

%

\section{Preliminaries}

Let $\frg_0$ be the Lie algebra of $G$. We always fix a Vogan diagram for $\frg_0$ as \cite[Appendix C]{Kn}. By doing this we have actually fixed a positive root system $\Delta^+(\frg, \frh_f)$. Restricting it to $\frt_f^*$, we get $(\Delta^+)^{(0)}(\frg, \frt_f)$ which contains a $\Delta^+(\frk, \frt_f)$. We fix this $\Delta^+(\frk, \frt_f)$ once for all. Let $\alpha_1, \dots, \alpha_l$ be the simple roots of $(\Delta^+)^{(0)}(\frg, \frt_f)$, and let $\zeta_1, \dots, \zeta_l$ be the corresponding fundamental weights.
There are $s$ ways of choosing  positive root systems for $\Delta(\frg, \frt_f)$ containing $\Delta^+(\frk, \frt_f)$, where $s=|W(\frg, \frt_f)|/|W(\frk, \frt_f)|$.
We enumerate them as
$$
(\Delta^+)^{(j)}(\frg, \frt_f)=\Delta^+(\frk, \frt_f) \cup  (\Delta^+)^{(j)}(\frp, \frt_f), \quad 0\leq j\leq s-1.
$$
Let us denote the half sum of roots of $(\Delta^+)^{(j)}(\frg, \frt_f)$ by $\rho^{(j)}$, and put $\rho_n^{(j)}:=\rho^{(j)} -\rho_c$. Then $\rho_n^{(j)}$ is the half sum of roots in $(\Delta^+)^{(j)}(\frp, \frt_f)$. Let $w^{(j)}$ be the unique element in $W(\frg, \frh_f)$ such that $w^{(j)} \rho^{(0)}=\rho^{(j)}$, and collect them as the set $W(\frg, \frt_f)^1$. Then $w^{(0)}=e$, and by a result of Kostant \cite{Ko}, the multiplication map induces a bijection from $W(\frg, \frt_f)^1 \times W(\frk, \frt_f)$ to $W(\frg, \frt_f)$.
For any $0\leq j \leq s-1$, $w^{(j)}\alpha_1, \dots, w^{(j)}\alpha_l$ are the simple roots of $(\Delta^+)^{(j)}(\frg, \frt_f)$, and $w^{(j)}\zeta_1, \dots, w^{(j)}\zeta_l$ are the corresponding fundamental weights.

Choose an irreducible module ${\rm Spin}_G$ of the Clifford algebra $C(\frp)$. Then as $\frk$-modules, we have the following decomposition
\begin{equation}\label{spin-module}
{\rm Spin}_G=2^{[\frac{l_0}{2}]} \bigoplus_{j=0}^{s-1} E_{\rho_n^{(j)}},
\end{equation}
where $l_0=\dim_{\bbC} \fra_f$. See Lemma 9.3.2 of \cite{Wa}.

Let us denote the dominant Weyl chamber for $(\Delta^+)^{(j)}(\frg, \frt_f)$ (resp., $\Delta^+(\frk, \frt_f)$) by $\caC^{(j)}$ (resp., $\caC$). Then
\begin{equation}\label{cone-deco}
\caC=\bigcup_{j=0}^{s-1}\caC^{(j)}.
\end{equation}
Note that $w^{(j)}\caC^{(0)}=\caC^{(j)}$ for $0\leq j\leq s-1$. The convex hull formed by the $W(\frk, \frt_f)$ orbits of the points $\rho_n^{(j)}$, $0\leq j\leq s-1$, is called the \emph{unitarily small convex hull}, and a $K$-type is said to be \emph{unitariy small} (\emph{u-small} for short henceforth) if its highest weight lies in this polyhedron. See Salamanca-Riba and Vogan \cite{SV}. Otherwise, we will call this $K$-type \emph{u-large}.

\subsection{Lambda norm and spin norm}\label{sec-lambda-spin-norm}
The following geometric way of describing the lambda norm in  \cite{V81} is due to Carmona \cite{Ca}. Consider an arbitrary $K$-type $E_{\mu}$. Choose $0\leq j\leq s-1$ such that $\mu + 2\rho_c$ is dominant for $(\Delta^+)^{(j)}(\frg, \frt_f)$.
Put
\begin{equation}\label{lambda-a-mu}
\lambda_a(\mu):=P(\mu+2 \rho_c -\rho^{(j)}).
\end{equation}
Here for any vector $\eta\in i\frt_{f, 0}^*$,  $P(\eta)$ stands for the projection of $\eta$ to the cone $\caC^{(j)}$. Namely, $P(\eta)$ is the unique point in $\caC^{(j)}$ which is closest to $\eta$. By Proposition 1.1(c) of \cite{SV}, the two vectors $\eta-P(\eta)$ and $P(\eta)$ are orthogonal to each other. Therefore,
\begin{equation}\label{proj-ortho}
\|\eta\|^2 -\|P(\eta)\|^2=\|\eta -P(\eta)\|^2.
\end{equation}
The lambda norm of $\mu$ is defined as
\begin{equation}\label{lambda-norm}
\|\mu\|_{\rm lambda}:=\|\lambda_a(\mu)\|.
\end{equation}
It turns out that $\|\mu\|_{\rm lambda}$ is independent of the choice of an allowable $j$ and is thus well-defined.

On the other hand, the spin norm introduced in \cite{D13} is
\begin{equation}\label{lambda-norm}
\|\mu\|_{\rm spin}:=\min_{0\leq j\leq s-1} \|\{\mu-\rho_n^{(j)}\}  + \rho_c \|.
\end{equation}
Here $\{\mu-\rho_n^{(j)}\}$ is the unique dominant weight in the $W(\frk, \frt_f)$ orbit of $\mu-\rho_n^{(j)}$. We emphasize that $\{\mu-\rho_n^{(j)}\}$, $0\leq j\leq s-1$, are precisely all the PRV components \cite{PRV} of the tensor product $E_{\mu}\otimes {\rm Spin}_G$ as $\frk$-modules. It has been shown in \cite{D13} that
\begin{equation}\label{lambda-norm}
\|\mu\|_{\rm spin}\geq \|\mu\|_{\rm lambda}.
\end{equation}

\subsection{Pathasarathy's Dirac operator inequality and a reformulation}
Let $Z_1, \dots, Z_m$ be an orthonormal basis of $\frp_0$ under the Killing form. Then
$$
D:=\sum_{i=1}^{m} Z_i\otimes Z_i \in U(\frg)\otimes C(\frp)
$$
is the Dirac operator introduced by Parthasarathy \cite{Pa72}. For any irreducible $(\frg, K)$-module $X$ with infinitesimal character \eqref{inf-char}, the Dirac operator $D$ acts on $X\otimes {\rm Spin}_G$. When $X$ is unitary, then $D^2\geq 0$ on $X\otimes {\rm Spin}_G$. Writing out $D^2$ carefully, one will get Parthasarathy's Dirac operator inequality \cite{Pa80}
\begin{equation}\label{Dirac-inequality}
\|\gamma + \rho_c\| \geq \|\Lambda\|,
\end{equation}
where $\gamma$ is the highest weight of any $\widetilde{K}$-type occurring in $X\otimes {\rm Spin}_G$. Here $\widetilde{K}$ is the spin covering group of $K$. This inequality is very effective in non-unitarity test: once a $\widetilde{K}$-type $\gamma$ in $X\otimes {\rm Spin}_G$ has been detected such that \eqref{Dirac-inequality} fails, then $\pi$ must be non-unitary. Finding such a $\gamma$ maybe skillful, yet the paper \cite{PRV} guarantees that it suffices to consider the PRV components of $X\otimes {\rm Spin}_G$.  Keeping this in mind, one sees that the inequality \eqref{Dirac-inequality} is equivalent to
\begin{equation}\label{Dirac-inequality-re}
\|\mu\|_{\rm spin} \geq \|\Lambda\|,
\end{equation}
where $\mu$ is the highest weight of any $K$-type of $\pi$.

\subsection{Non-decreasable $K$-types should be u-small}

For any $0\leq j \leq s-1$, collect all the $K$-types $\mu$ such that
$\mu+2\rho_c\in \caC^{(j)}$ as $\Omega(j)$. We call a $K$-type $\mu\in\Omega(j)$ \emph{decreasable} if $\mu-w^{(j)}\zeta_p\in\caC^{(j)}$ for certain $1\leq p\leq l$. Otherwise, $\mu$ is called \emph{non-decreasable}.

\begin{conj}\label{conj-non-dec}
Assume that the $K$-type $\mu\in \Omega(j)$ is non-decreasable. Then $\mu$ should be u-small.
\end{conj}

We will give explicit examples of non-decreasable $K$-types in Sections \ref{sec-real-exc} and \ref{sec-HS}.
The above conjecture has been verified for all exceptional groups that we shall address.

\section{Sharpening the Helgason-Johnson bound}\label{sec-framework}

Learning from Vogan's lectures \cite{V20-1,V20-2}, this section aims to give a way to improve the Helgason-Johnson bound for infinite-dimensional unitary representations.
Indeed, let $\pi$ be such a representation with a lowest $K$-type $E_{\mu}$, and let its infinitesimal character be given as in \eqref{inf-char}. Now by \eqref{Dirac-inequality-re},
$$
\|\Lambda\|^2=\|\lambda_a(\mu)\|^2+\|\nu\|^2\leq \|\mu+m \beta\|^2_{\rm spin}, \quad \forall m\in \bbZ_{\geq 0}.
$$
Therefore,
$$
\|\nu\|^2\leq \|\mu+m \beta\|^2_{\rm spin}- \|\lambda_a(\mu)\|^2, \quad \forall m\in \bbZ_{\geq 0}.
$$
In other words,
\begin{equation}\label{nu-bound}
\|\nu\|^2\leq {\rm MP}(\mu)- \|\mu\|_{\rm lambda}^2.
\end{equation}
To get a uniform bound of $\|\nu\|^2$ for all the infinite-dimensional irreducible unitary representations, it suffices to obtain the \textbf{maximum} of the RHS of \eqref{nu-bound} over all the $K$-types $\mu$.

Now let us explain the idea of handling the RHS of \eqref{nu-bound}. Theorem 1.1 of \cite{D16} and Theorem C of \cite{D17} suggest that when $\mu$ is u-large, we should have ${\rm MP(\mu)}=\|\mu\|_{\rm spin}^2$. Thus we should firstly understand the behavior of $\|\mu\|_{\rm spin}^2 - \|\mu\|_{\rm lambda}^2$.

Take any $\mu\in \Omega(j)$. We have that
\begin{align*}
\|\mu\|_{\rm spin}^2 - \|\mu\|_{\rm lambda}^2
&=\|\mu\|_{\rm spin}^2-\|P(\mu+2\rho_c-\rho^{(j)})\| \\
&\leq \|\{\mu-\rho_n^{(j)}\}+\rho_c\|^2-\|P(\mu+\rho_c-\rho_n^{(j)})\| \\
&\leq {\rm I}(\mu) + {\rm II }(\mu),
\end{align*}
where
\begin{equation}\label{term-I}
{\rm I}(\mu)=\|\{\mu-\rho_n^{(j)}\}+\rho_c\|^2 -\|\mu-\rho_n^{(j)}+\rho_c\|^2=2\langle \rho_c,
\{\mu-\rho_n^{(j)}\} - (\mu-\rho_n^{(j)})\rangle
\end{equation}
and
\begin{equation}\label{term-II}
{\rm II}(\mu)=\|\mu+\rho_c -\rho_n^{(j)}\|^2-\|P(\mu+\rho_c-\rho_n^{(j)})\|^2.
\end{equation}
Here recall that $P(\mu+\rho_c-\rho_n^{(j)})$ is the projection of $\mu+\rho_c-\rho_n^{(j)}$ to the cone $\caC^{(j)}$. By \eqref{proj-ortho},
\begin{equation}\label{term-II-equiv}
{\rm II}(\mu)=\|\mu+\rho_c -\rho_n^{(j)}-P(\mu+\rho_c-\rho_n^{(j)})\|^2.
\end{equation}
It follows that both ${\rm I}(\mu)$ and ${\rm II}(\mu)$ are non-negative.

Take any vector $\delta\in\caC^{(j)}$, Lemma \ref{lemma-neg-comp} will guarantee that ${\rm I}(\mu+\delta)\leq {\rm I}(\mu)$ since $\delta$ must be dominant for $\Delta^+(\frk, \frt_f)$. On the other hand,
\begin{align*}
{\rm II}(\mu+\delta) &= \|\mu+\delta+\rho_c -\rho_n^{(j)}- P(\mu+\delta+\rho_c -\rho_n^{(j)})\|^2 \\
                     &\leq \|\mu+\delta+\rho_c -\rho_n^{(j)}-(\delta + P(\mu+\rho_c-\rho_n^{(j)}))\|^2 \\
                     &= \|\mu+\rho_c-\rho_n^{(j)}-P(\mu+\rho_c-\rho_n^{(j)})\|^2 \\
                     &={\rm II}(\mu),
\end{align*}
where the second step holds since $\delta + P(\mu+\rho_c-\rho_n^{(j)})\in\caC^{(j)}$.

Now we can state the following way of sharpening the Helgason-Johnson bound:

(a) For any $0\leq j \leq s-1$, select finitely many u-small $K$-types $\mu_1, \dots, \mu_{m(j)}$ in $\Omega(j)$ so that for any u-large $\mu\in \Omega(j)$, there exists $1\leq i\leq m(j)$ such that $\mu-\mu_i\in \caC^{(j)}$.
    Then as deduced above,
    $$
    {\rm I}(\mu) +  {\rm II}(\mu) \leq  {\rm I}(\mu_i) +  {\rm II}(\mu_i).
    $$
    Let
    \begin{equation}\label{bound-Aj}
    A_j:=\max\{{\rm I}(\mu_i) +  {\rm II}(\mu_i)\mid 1\leq i\leq m(j)\}.
    \end{equation}
    Therefore, for any u-large $\mu\in\Omega(j)$, we have that $\|\mu\|_{\rm spin}^2 -\|\mu\|_{\rm lambda}^2 \leq A_j$.

(b) We directly compute
\begin{equation}\label{bound-B}
B=\max_{\mu}\left\{\min_m \|\mu+m\beta\|^2_{\rm spin} - \|\mu\|^2_{\rm lambda}\mid m\in\bbZ_{\geq 0} \mbox{ such that } \mu+m\beta \mbox{ is u-small} \right\},
\end{equation}
where $\mu$ runs over all the u-small $K$-types in $\Omega(j)$ other than those $\mu_1, \dots, \mu_{m(j)}$ in step (a).

To sum up, $\|\nu\|^2$ is upper bounded by $\max\{A_0, A_1, \dots, A_{s-1}, B\}$.


\subsection{Technical lemmas for ${\rm I}(\mu)$}
This subsection aims to prepare some technical lemmas for analyzing the term ${\rm I}(\mu)$ defined in \eqref{term-I}. Let $\gamma_1, \dots, \gamma_l$ be the simple roots of $\Delta^+(\frk, \frt_f)$, and let $\varpi_1, \dots, \varpi_l$ be the corresponding fundamental weights.

Let $\eta=\sum_{i=1}^l \eta_i\varpi_i$ be
an arbitrary weight. We utilize the following \emph{negative index
algorithm} to conjugate $\eta$ to $\caC$: if the weight is not yet
dominant, select an index $i$ such that $\eta_i<0$, then apply the simple reflection $s_{\gamma_i}$.  We record the corresponding sequence of
simple reflections, and collect all of them as $S(\lambda)$.

\begin{lemma}\label{lemma-reduced} \emph{(Theorem 4.3.1(iv) of Bj\"orner and Brenti \cite{BB}, Proposition 4.1 of Eriksson \cite{Er})}
Let $\eta$ be any weight. Let $s_{\beta_1}, \dots, s_{\beta_n}, \dots$
be any sequence in $S(\eta)$, then the word $s_{\beta_n} \cdots
s_{\beta_1}$ is reduced. Here each $\beta_i$ is a simple root of $\Delta^+(\frk, \frt_f)$. In particular, each sequence in $S(\eta)$ is of finite length.
\end{lemma}

\begin{lemma}\label{lemma-neg-ext} \emph{(Theorem 4.5 of Eriksson \cite{Er})}
Assume that $\eta_1$ and $\eta_2$ are two weights such that $\eta_1-\eta_2$ is dominant for $\Delta^+(\frk, \frt_f)$. Then every sequence $s_{\beta_1}, \dots, s_{\beta_n}$ in $S(\eta_1)$ can be extended to a sequence $s_{\beta_1}, \dots, s_{\beta_n}, \dots, s_{\beta_N}$ in $S(\eta_2)$. Here $N\geq n$.
\end{lemma}
\begin{proof}
Note that for any $1\leq k\leq n$, we have
$$
\langle \eta_1 -\eta_2,  s_{\beta_1}s_{\beta_2}\cdots s_{\beta_{k-1}}(\beta_k)\rangle\geq 0
$$
since $\eta_1-\eta_2$ is dominant for $\Delta^+(\frk, \frt_f)$ and $s_{\beta_1}s_{\beta_2}\cdots s_{\beta_{k-1}}(\beta_k)\in \Delta^+(\frk, \frt_f)$. Therefore,
$$
\langle s_{\beta_{k-1}}\cdots s_{\beta_2} s_{\beta_{1}}(\eta_1 -\eta_2),  \beta_k\rangle\geq 0.
$$
Thus for any $1\leq k\leq n$,
$$
\langle s_{\beta_{k-1}}\cdots s_{\beta_2} s_{\beta_{1}}(\eta_2),  \beta_k^\vee\rangle\leq \langle s_{\beta_{k-1}}\cdots s_{\beta_2} s_{\beta_{1}}(\eta_1),  \beta_k^\vee\rangle<0.
$$
Now the negative index algorithm finishes the proof.
\end{proof}

\begin{lemma}\label{lemma-neg-comp}
Assume that $\eta_1$ and $\eta_2$ are two weights such that $\eta_1-\eta_2$ is dominant for $\Delta^+(\frk, \frt_f)$. Then $\langle \rho_c, \{\eta_1\} -\eta_1 \rangle\leq \langle \rho_c, \{\eta_2\} -\eta_2 \rangle$.
\end{lemma}
\begin{proof}
Take a sequence $s_{\beta_1}, \dots, s_{\beta_n}$ in $S(\eta_1)$. By Lemma \ref{lemma-neg-ext}, it extends to a sequence $s_{\beta_1}, \dots, s_{\beta_n}, \dots, s_{\beta_N}$ in $S(\eta_2)$. Note that $s_{\beta_n}\cdots s_{\beta_1}$ and $s_{\beta_N}\cdots s_{\beta_n}\cdots s_{\beta_1}$  are reduced expressions due to Lemma \ref{lemma-reduced}.

By Lemma 5.5 of \cite{DH},
\begin{align*}
\{\eta_1\} -\eta_1&=s_{\beta_n}\cdots s_{\beta_1}(\eta_1)-\eta_1\\
                  &=\sum_{k=1}^n \left( s_{\beta_n}\cdots s_{\beta_k}(\eta_1) -s_{\beta_n}\cdots s_{\beta_{k+1}}(\eta_1)  \right)\\
                  &=\sum_{k=1}^n  s_{\beta_n}\cdots s_{\beta_{k+1}}(\eta_1-s_{\beta_{k}}(\eta_1))  \\
                  &=\sum_{k=1}^n \langle -\eta_1, \beta_k^\vee \rangle   s_{\beta_n}\cdots s_{\beta_{k+1}}(\beta_{k}).
\end{align*}
Therefore,
\begin{equation}\label{eta-1-rhoc}
\langle \rho_c, \{\eta_1\} -\eta_1 \rangle=\sum_{k=1}^n
\langle -\eta_1, \beta_k^\vee \rangle  \langle \rho_c,  s_{\beta_n}\cdots s_{\beta_{k+1}}(\beta_{k}) \rangle.
\end{equation}

On the other hand,
\begin{align*}
\{\eta_2\} -\eta_2&=
\sum_{k=1}^n \left( s_{\beta_n}\cdots s_{\beta_k}(\eta_2) -s_{\beta_n}\cdots s_{\beta_{k+1}}(\eta_2)  \right)\\
                  &+ \sum_{k=n+1}^N \left(
                  s_{\beta_k}\cdots s_{\beta_1}(\eta_2) -s_{\beta_{k-1}}\cdots s_{\beta_{1}}(\eta_2)\right)\\
                  &=\sum_{k=1}^n \langle -\eta_2, \beta_k^\vee \rangle   s_{\beta_n}\cdots s_{\beta_{k+1}}(\beta_{k})-\sum_{k=n+1}^N \langle s_{\beta_{k-1}}\cdots s_{\beta_1}(\eta_2), \beta_k^\vee \rangle \beta_k
\end{align*}
Therefore,
\begin{align*}
\langle \rho_c, (\{\eta_2\} -\eta_2) -(\{\eta_1\} -\eta_1) \rangle&=\sum_{k=1}^n
\langle \eta_1-\eta_2, \beta_k^\vee \rangle  \langle \rho_c,  s_{\beta_n}\cdots s_{\beta_{k+1}}(\beta_{k}) \rangle \\
 &-\sum_{k=n+1}^N \langle s_{\beta_{k-1}}\cdots s_{\beta_1}(\eta_2), \beta_k^\vee \rangle \langle \rho_c, \beta_k\rangle.
\end{align*}
Since the sequence $s_{\beta_1}, \dots, s_{\beta_n}, \dots, s_{\beta_N}$ belongs to $S(\eta_2)$, we have that
$$
\langle s_{\beta_{k-1}}\cdots s_{\beta_1}(\eta_2), \beta_k^\vee \rangle < 0, \quad \forall n+1\leq k\leq N
$$
Now the desired conclusion follows.
\end{proof}

\section{Complex exceptional Lie groups}
This section aims to prove Theorem \ref{thm-complex-exceptional}.
Assume that $G$ is complex, let $\alpha_1, \dots, \alpha_l$ be the simple roots for the fixed $\Delta^+(\frk, \frt_f)$. Let $\varpi_1, \dots, \varpi_l$ be the corresponding fundamental weights. The highest weight $\mu$ of a $\frk$-type is a non-negative integer combination of $\varpi_1, \dots, \varpi_l$. For convenience, we use $[a_1, a_2, \dots, a_l]$ to stand  for the vector $\sum_{i=1}^{l} a_i \varpi_i$.

Note that $s=1$, $\rho^{(0)}=2\rho_c$ and $\rho_n^{(0)}=\rho_c$. Collect all the u-small $\frk$-types as $\Omega_{\rm us}(0)$. Let $B$ be the maximum of the expression ${\rm MP}(\mu)-\|\mu\|^2_{\rm lambda}$  over the set $\Omega_{\rm us}(0)$.
Define
\begin{equation}\label{bdry-Omega-complex}
\partial\Omega_{\rm us}(0)=\{\mu\in\Omega_{\rm us}(0)\mid \exists\, 1\leq i \leq l \mbox{ such that } \mu + \varpi_i \mbox{ is not u-small}\}.
\end{equation}
Let $A_0$ be the maximum of the expression $\|\mu\|^2_{\rm spin}-\|\mu\|^2_{\rm lambda}$ over the set $\partial\Omega_{\rm us}(0)$. Note that for any u-large $\frk$-type $\mu$, there exists $\mu_1\in\partial\Omega_{\rm us}(0)$ such that $\mu-\mu_1\in\caC$. Indeed, we can subtract some $\varpi_i$ from $\mu$ step by step, but stay within the cone $\caC$. Then by Conjecture \ref{conj-non-dec} (which is easily seen to be true for complex groups), eventually we will come down to the u-small convex hull. The first $\frk$-type that we meet in the u-small convex hull can be chosen as the desired $\mu_1$. Therefore, due to the discussion around \eqref{bound-B}, $\max\{A_0, B\}$ is an upper bound of $\|\nu\|^2$ in Theorem \ref{thm-complex-exceptional}.

We will handle the complex exceptional groups one by one. It turns out that $B$ is attained at the trivial $\frk$-type in each case.

\subsection{Complex $G_2$}

\begin{figure}[H]
\centering
\scalebox{0.5}{\includegraphics{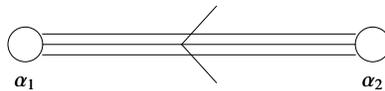}}
\caption{The Dynkin diagram for $G_2$}
\label{Fig-G2-Dynkin}
\end{figure}

Let $\alpha_1=(1, -1, 0)$ and $\alpha_2=(-2, 1, 1)$ be the two simple roots. The Dynkin diagram is presented in Fig.~\ref{Fig-G2-Dynkin}.

\begin{figure}[H]
\centering
\scalebox{0.80}{\includegraphics{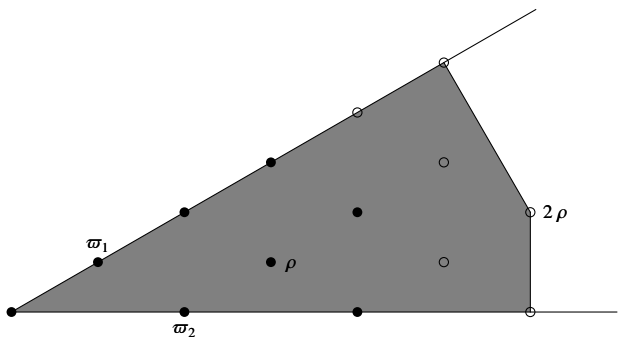}}
\caption{The complex $G_2$ case}
\label{Fig-G2-USmall}
\end{figure}

In this case, $\Omega_{\rm us}(0)$ has fourteen $\frk$-types and
$$
\partial\Omega_{\rm us}(0)=\{[0, 3], [1, 2], [2, 2], [3, 1], [4, 0], [5, 0]\}.
$$
In Fig.~\ref{Fig-G2-USmall}, we use circles to denote the $\frk$-types in $\partial\Omega_{\rm us}(0)$, while the other u-small $\frk$-types are denoted by dots.

We compute that $A_0=6$ and $B=26$, which is attained at the trivial $\frk$-type.
 Thus the complex $G_2$ case in Theorem \ref{thm-complex-exceptional} follows.

\subsection{Complex $F_4$}

\begin{figure}[H]
\centering
\scalebox{0.5}{\includegraphics{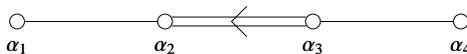}}
\caption{The Dynkin diagram for $F_4$}
\label{Fig-F4-Dynkin}
\end{figure}

Let $\alpha_1=(\frac{1}{2}, -\frac{1}{2}, -\frac{1}{2}, -\frac{1}{2})$, $\alpha_2=(0, 0, 0, 1)$, $\alpha_3=(0, 0, 1, -1)$ and $\alpha_4=(0, 1, -1, 0)$ be the simple roots. The Dynkin diagram is presented in Fig.~\ref{Fig-F4-Dynkin}.

We compute that $|\Omega_{\rm us}(0)|=451$ and $|\partial\Omega_{\rm us}(0)|=214$. Moreover, $A_0=35$ and $B=78$, which is attained at the trivial $\frk$-type.
The complex $F_4$ case in Theorem \ref{thm-complex-exceptional} follows.

\subsection{Complex $E_6$}
\begin{figure}[H]
\centering
\scalebox{0.6}{\includegraphics{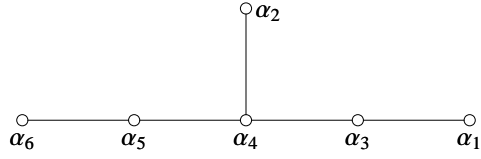}}
\caption{The Dynkin diagram for $E_6$}
\label{Fig-E6-Dynkin}
\end{figure}

We label the Dynkin diagram in Fig.~\ref{Fig-E6-Dynkin}.

We compute that $|\Omega_{\rm us}(0)|=13007$ and $|\partial\Omega_{\rm us}(0)|=6249$. Moreover, $A_0=120$ and $B=170$, which is attained at the trivial $\frk$-type.
The complex $E_6$ case in Theorem \ref{thm-complex-exceptional} follows.

\subsection{Complex $E_7$}
\begin{figure}[H]
\centering
\scalebox{0.6}{\includegraphics{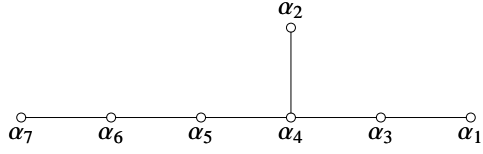}}
\caption{The Dynkin diagram for $E_7$}
\label{Fig-E7-Dynkin}
\end{figure}

We label the Dynkin diagram in Fig.~\ref{Fig-E7-Dynkin}.

We compute that $|\Omega_{\rm us}(0)|=105951$ and $|\partial\Omega_{\rm us}(0)|=52630$. Moreover, $A_0=312$ and $B=464$, which is attained at the trivial $\frk$-type. The complex $E_7$ case in Theorem \ref{thm-complex-exceptional} follows.

\subsection{Complex $E_8$}
\begin{figure}[H]
\centering
\scalebox{0.65}{\includegraphics{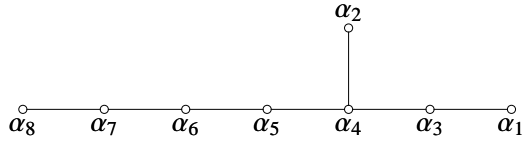}}
\caption{The Dynkin diagram for $E_8$}
\label{Fig-E8-Dynkin}
\end{figure}

We label the Dynkin diagram in Fig.~\ref{Fig-E8-Dynkin}.

We compute that $|\Omega_{\rm us}(0)|=950077$ and $|\partial\Omega_{\rm us}(0)|=486023$. Moreover, $A_0=800$ and $B=1522$, which is attained at the trivial $\frk$-type. The complex $E_8$ case in Theorem \ref{thm-complex-exceptional} follows.

\section{Real exceptional Lie groups}\label{sec-real-exc}

This section aims to prove Theorem \ref{thm-real-exceptional}. Recall that $\alpha_1, \dots, \alpha_l$ are the simple roots of $(\Delta^+)^{(0)}(\frg, \frt_f)$, with $\zeta_1, \dots, \zeta_l$ being the corresponding fundamental weights.
For convenience, we will use $s_i$ to denote the simple reflection $s_{\alpha_i}$, and use $[[a_1, \dots, a_l]]$ to denote the vector $a_1\zeta_1+\cdots+ a_l\zeta_l$.

Since $\frk$ is assumed to have no center, $\Delta^+(\frk, \frt_f)$ has $l$ simple roots as well. We denote them by $\gamma_1, \dots, \gamma_l$, and let the corresponding fundamental weights be $\varpi_1, \dots, \varpi_l$. The highest weight of any $\frk$-type is a non-negative integer combination of $\varpi_1, \dots, \varpi_l$. For convenience, we will write $[b_1, \dots, b_l]$ to stand for the vector $b_1\varpi_1+\cdots+ b_l\varpi_l$.

Fix any $0\leq j\leq s-1$.
Collect all the $\frk$-types $\mu$ such that $\mu+2\rho_c\in\caC^{(j)}$ as $\Omega(j)$.
Let $\Omega_{\rm us}(j)$ consist of the u-small members of $\Omega(j)$. Define
\begin{equation}\label{bdry-Omega-real}
\partial\Omega_{\rm us}(j)=\{\mu\in\Omega_{\rm us}(j)\mid \exists\, 1\leq i \leq l \mbox{ such that } \mu + w^{(j)}\zeta_i \mbox{ is not u-small}\}.
\end{equation}
Let $A_j$ be the maximum of  $\|\{\mu-\rho_n^{(j)}\}+\rho_c\|^2-\|\mu\|^2_{\rm lambda}$ over the set $\partial\Omega_{\rm us}(j)$.

We \textbf{claim} that for any u-large $\frk$-type $\mu\in \Omega(j)$, there exists $\mu_1\in \partial\Omega_{\rm us}(j)$ such that $\mu-\mu_1\in\caC^{(j)}$.
Indeed, it follows from  Conjecture \ref{conj-non-dec} (once verified) that $\mu\in\Omega(j)$ is decreasable. Namely, we can subtract some $w^{(j)}\zeta_i$ from $\mu$ so that $\mu-w^{(j)}\zeta_i\in\Omega(j)$. Doing this step by step, eventually we will come down to a non-decreasable $\frk$-type in $\Omega(j)$ which must be u-small. The first $\frk$-type that we meet in the u-small convex hull during this going down process can be chosen as the desired $\mu_1$. Thus the claim holds.

Let $B$ be the maximum of the expression
$$
\min_m \{\|\mu+m\beta\|^2_{\rm spin}-\|\mu\|^2_{\rm lambda}\mid m\in\bbZ_{\geq 0} \mbox{ such that } \mu+m\beta \mbox{ is u-small}\}
$$ with $\mu$ running over $\Omega_{\rm us}(j)\setminus\partial \Omega_{\rm us}(j)$.
Due to the discussion around \eqref{bound-B}, $\max\{A_0, \dots, A_{s-1}, B\}$ is an upper bound of $\|\nu\|^2$ in Theorem \ref{thm-real-exceptional}.

In the following, we will handle the real exceptional Lie groups one by one. As mentioned in the introduction, the value $B$ is attained at the trivial $K$-type if and only if $G/K$ is \emph{not} quaternionic in the sense of \cite{Kn}.

\subsection{$G_{2(2)}$}

The Vogan diagram of this simple Lie algebra is obtained from Fig.~\ref{Fig-G2-Dynkin} by painting the simple root $\alpha_2$. Then $\gamma_1=\alpha_1$ and $\gamma_2=3\alpha_1+ 2 \alpha_2$  are the simple roots for $\Delta^+(\frk, \frt_f)$, which has type $A_1 \times A_1$. Let $\mu=[a, b]$ be a $\frk$-type.
Let us enumerate the set $W(\frg, \frt_f)^1$ in order as $w^{(0)}=e$, $w^{(1)}=s_2$ and $w^{(2)}=s_2s_1$. Then
$$\rho_n^{(0)}=[0, 2], \quad \rho_n^{(1)}=\beta=[3, 1], \quad \rho_n^{(2)}=[4, 0].
$$
Moreover, $\mu+2\rho_c$ has coordinates
$$
(a+2, \frac{b}{2} -\frac{a}{2}), \quad (-\frac{a}{2}+\frac{3b}{2}+2, \frac{a}{2} -\frac{b}{2}), \quad (\frac{a}{2}-\frac{3b}{2}-2, b+2)
$$
with respect to the basis $\{w^{(j)}\zeta_1, w^{(j)}\zeta_2\}$ for $j=0, 1, 2$, respectively.

We move on according to the following three cases:
\begin{itemize}
 \item[$\bullet$]  $\mu\in\Omega(0)$ if and only if $b\geq a$. Note that $\zeta_1=[1, 1]$ and $\zeta_2=[0, 2]$. We compute that
     $$
     \partial\Omega_{\rm us}(0)=\{[0, 3], [0, 4], [1, 2], [1, 3], [2, 2], [2, 3], [3, 3]\}
     $$
     and that $A_0=\frac{3}{2}$.

\item[$\bullet$]  $\mu\in\Omega(1)$ if and only if $b\leq a \leq 3b+4$.

Note that $w^{(1)}\zeta_1=[1, 1]$ and $w^{(1)}\zeta_2=[3, 1]$. We compute that
     $$
     \partial\Omega_{\rm us}(1)=\{[4, 1], [5, 1], [6, 1], [7, 1], [2, 2], [3, 2], [4, 2], [5, 2], [6, 2], [3, 3]\}
     $$
     and that $A_1=\frac{7}{2}$.

\item[$\bullet$]  $\mu\in\Omega(2)$ if and only if $a\geq 3b+4$.

Note that $w^{(2)}\zeta_1=[2, 0]$ and $w^{(2)}\zeta_2=[3, 1]$. We compute that
     $$
     \partial\Omega(2)_{\rm us}=\{[5, 0], [6, 0], [7, 0], [8, 0],  [7, 1]\}
     $$
     and that $A_2=\frac{1}{2}$.
\end{itemize}
Moreover, we compute that $B=\frac{51}{8}$, which happens at the $\frk$-type $[0, 1]$. The $G_{2(2)}$ case in Theorem \ref{thm-real-exceptional} follows.

\begin{figure}[H]
\centering
\scalebox{0.70}{\includegraphics{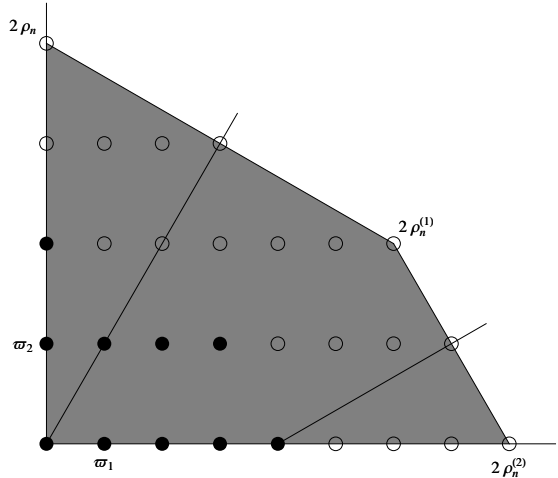}}
\caption{The $G_{2(2)}$ case}
\label{Fig-G-USmall}
\end{figure}

In Fig.~\ref{Fig-G-USmall}, the $\frk$-types in $\partial\Omega_{\rm us}(i)$ for $0\leq i\leq 2$ are denoted by circles, while the other u-small $\frk$-types are denoted by dots.

\begin{example}\label{exam-G}
The non-decreasable $\frk$-types in $\Omega(j)$ are
$$
j=0: [0, 0], [0, 1]; \quad
j=1: [0, 0], [1, 0], [2, 0], [3, 0], [4, 0]; \quad
j=2: [4, 0], [5, 0].
$$
They are all u-small. \hfill\qed
\end{example}

\subsection{${\rm FI}=F_{4(4)}$}
The Vogan diagram of this simple Lie algebra is obtained from Fig.~\ref{Fig-F4-Dynkin} by painting the simple root $\alpha_4$. Then $\gamma_i=\alpha_i$ for $1\leq i\leq 3$ and $\gamma_4:=2\alpha_1+4\alpha_2+3\alpha_3+2\alpha_4$ are the simple roots for $\Delta^+(\frk, \frt_f)$, which has type $C_3 \times A_1$. Note that $\beta=[0, 0, 1, 1]$.

\begin{table}[H]
\centering
\caption{The $F_{4(4)}$ case}
\begin{tabular}{l|r|c|c|c|r}
$j$ &   $w^{(j)}$ &$\rho_n^{(j)}$  & $\#\Omega_{\rm us}(j)$ & $\#\partial\Omega_{\rm us}(j)$ & $A_j$\\
\hline
$0$ & $e$   &$[0, 0, 0, 7]$ & $12$ & $11$ & $0.5$\\
$1$ & $s_4$ &$[0, 0, 1, 6]$ & $68$ & $53$ & $8.5$\\
$2$ & $s_4s_3$ &$[0, 2, 0, 5]$ & $116$ & $84$ & $9.5$\\
$3$ & $s_4s_3s_2$     &$[1, 2, 0, 4]$ & $193$ & $132$ & $9.5$\\
$4$ & $s_4s_3s_2s_1$  &$[0, 3, 0, 3]$ & $105$ & $74$ & $7$\\
$5$ & $s_4s_3s_2s_3$  &$[3, 0, 1, 3]$ & $166$ & $114$ & $9$\\
$6$ & $s_4s_3s_2s_1s_3$ &$[2, 1, 1, 2]$ & $339$ & $215$ & $9$\\
$7$ & $s_4s_3s_2s_3s_4$ &$[5, 0, 0, 2]$ & $41$ & $34$ & $0.75$\\
$8$ & $s_4s_3s_2s_1s_3s_2$ &$[2, 0, 2, 1]$ & $215$ & $142$ & $9.75$\\
$9$ & $s_4s_3s_2s_1s_3s_4$ &$[4, 1, 0, 1]$ & $130$ & $94$ & $6.25$\\
$10$ & $s_4s_3s_2s_1s_3s_2s_3$ &$[0, 0, 3, 0]$ & $43$ & $34$ & $2.5$\\
$11$ & $s_4s_3s_2s_1s_3s_2s_4$  &$[4, 0, 1, 0]$ & $87$ & $63$ & $6.25$
\end{tabular}
\label{table-FI}
\end{table}

We summarize the information in Table \ref{table-FI}. Moreover, $B=\frac{39}{2}$ and it is attained at the $\frk$-type $\mu=[0, 0, 0, 1]$. Indeed, $\mu$ belongs to $\Omega(6)$,
and $\mu+2\rho_c-\rho^{(6)}$ has coordinates
$(0, 0, -\frac{1}{2}, \frac{1}{2})$
w.r.t. the basis $\{w^{(6)}\zeta_1, \dots, w^{(6)}\zeta_4\}$. Therefore, projecting
$\mu+2\rho_c-\rho^{(6)}$ to $\caC^{(6)}$ is equivalent to projecting $[[0, 0, -\frac{1}{2}, \frac{1}{2}]]$ to $\caC^{(0)}$, which produces the zero vector. Therefore, $\|\mu\|_{\rm lambda}=0$. On the other hand, ${\rm MP}(\mu)=\|\mu+\beta\|_{\rm spin}=\frac{39}{2}$.  The $F_{4(4)}$ case in Theorem \ref{thm-real-exceptional} follows.

\begin{example}\label{exam-FI}
Take $j=6$ and consider the $\frk$-type $\mu=[a, b, c, d]$. Then $\mu+2\rho_c$ has coordinates
$$
(b+c-d+2, d-c, a/2+c/2-d/2+1, -a/2+c/2+d/2+1)
$$
w.r.t. the basis $\{w^{(j)}\zeta_1, \dots, w^{(j)}\zeta_4\}$. Therefore,
$\mu\in\Omega(j)$ if and only if
$$
d \leq b+c+2, \quad c\leq d, \quad d \leq a+c+2, \quad a\leq c+d+2.
$$
On the other hand,
$$
w^{(j)}\zeta_1=[0, 1, 0, 0], \quad w^{(j)}\zeta_2=[1, 1, 0, 1], \quad
w^{(j)}\zeta_3=[2, 0, 1, 1], \quad w^{(j)}\zeta_4=[0, 0, 1, 1].
$$
We compute that the non-decreasable $\frk$-types in $\Omega(j)$ are
\begin{align*}
[0,0,0,0], [0,0,0,1], [0,0,0,2], [1,0,0,0], [1,0,0,1], [1,0,0,2], \\
[2,0,0,0], [2,0,0,1], [2,0,0,2], [3,0,0,1], [3,0,0,2], [4,0,0,2].
\end{align*}
They are all u-small.  \hfill\qed
\end{example}

\subsection{${\rm FII}=F_{4(-20)}$}
The Vogan diagram of this simple Lie algebra is obtained from Fig.~\ref{Fig-F4-Dynkin} by painting the simple root $\alpha_1$. Then $\gamma_1=(1, -1, 0, 0)$, $\gamma_2=\alpha_4$, $\gamma_3=\alpha_3$ and $\gamma_4=\alpha_2$  are the simple roots for $\Delta^+(\frk, \frt_f)$, which is of type $B_4$. Let $\mu=[a, b, c, d]$ be a $\frk$-type. Note that $\beta=[0, 0, 0, 1]$.

\begin{table}[H]
\centering
\caption{The $F_{4(-20)}$ case}
\begin{tabular}{l|r|c|c|c|r}
$j$ &   $w^{(j)}$ &$\rho_n^{(j)}$  & $\#\Omega_{\rm us}(j)$ & $\#\partial\Omega_{\rm us}(j)$ & $A_j$\\
\hline
$0$ & $e$   &$[2, 0, 0, 0]$ & $5$ & $5$ & $0.25$\\
$1$ & $s_1$ &$[1, 0, 0, 1]$ & $21$ & $18$ & $11.25$\\
$2$ & $s_1s_2$ &$[0, 0, 1, 0]$ & $14$ & $12$ & $10.25$
\end{tabular}
\label{table-FII}
\end{table}

We summarize the information in Table \ref{table-FII}. Moreover, $B=\frac{81}{4}$, which is attained at the trivial $\frk$-type. The $F_{4(-20)}$ case in Theorem \ref{thm-real-exceptional} follows.


\subsection{${\rm EI}=E_{6(6)}$}
\begin{figure}[H]
\centering
\scalebox{0.60}{\includegraphics{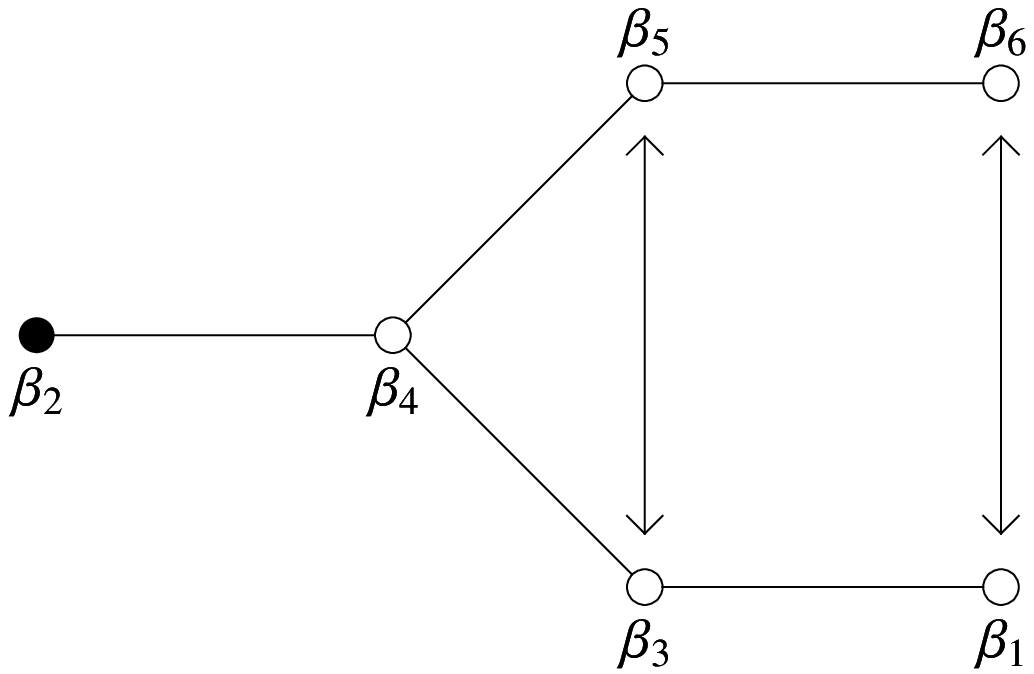}}
\caption{The Vogan diagram for EI}
\label{Fig-EI-Vogan}
\end{figure}

The Vogan diagram of this simple Lie algebra is presented in Fig.~\ref{Fig-EI-Vogan}.
The positive root system $(\Delta^+)^{(0)}(\frg, \frt_f)$ is of type $F_4$, and has simple roots $\alpha_1=\frac{1}{2}(\beta_1+\beta_6)$, $\alpha_2=\frac{1}{2}(\beta_3+\beta_5)$, $\alpha_3=\beta_4$, $\alpha_4=\beta_2$. Here $\alpha_1$ and $\alpha_2$ are short, while $\alpha_3$ and $\alpha_4$ are long.  On the other hand, the positive root system $\Delta^+(\frk, \frt_f)$ is of type $C_4$, and has simple roots $\gamma_1=\alpha_2+\alpha_3+\alpha_4$, $\gamma_2=\alpha_1$, $\gamma_3=\alpha_2$, $\gamma_4=\alpha_3$. Here $\gamma_4$ is long. Note that $\beta=[0, 0, 0, 1]$.

We summarize the information in Table \ref{table-EI}. Moreover, $B=42$, which is attained at the trivial $\frk$-type. The $E_{6(6)}$ case in Theorem \ref{thm-real-exceptional} follows.

\begin{table}[H]
\centering
\caption{The $E_{6(6)}$ case}
\begin{tabular}{l|r|c|c|c|r}
$j$ &   $w^{(j)}$ &$\rho_n^{(j)}$  & $\#\Omega_{\rm us}(j)$ & $\#\partial\Omega_{\rm us}(j)$ & $A_j$\\
\hline
$0$ & $e$   &$[5, 1, 1, 0]$ & $192$ & $124$ & $16$\\
$1$ & $s_4$ &$[3, 1, 1, 1]$ & $539$ & $295$ & $18.5$\\
$2$ & $s_4s_3$ &$[1, 1, 3, 0]$ & $354$ & $202$ & $16.5$
\end{tabular}
\label{table-EI}
\end{table}

\subsection{${\rm EII}=E_{6(2)}$}
The Vogan diagram of this simple Lie algebra is obtained by painting the simple root $\alpha_2$ in Fig.~\ref{Fig-E6-Dynkin}. Then $\Delta^+(\frk, \frt_f)$ is of type $A_5\times A_1$, and has the following simple roots
$$
\gamma_i=\alpha_{7-i}, \quad 1\leq i\leq 4;\quad \gamma_5=\alpha_1; \quad \gamma_6=\alpha_1+2\alpha_2+2\alpha_3+3\alpha_4+2\alpha_5+\alpha_6.
$$
We have that $|W(\frg, \frt_f)^1|=36$ and $\beta=[0, 0, 1, 0, 0, 1]$.

\begin{table}[H]
\centering
\caption{The $E_{6(2)}$ case}
\begin{tabular}{c|c|c|c|c|c|c|c|c|c}
$\rho_n^{(j)}$ &  $\#\Omega_{\rm us}(j)$ & $\#\partial\Omega_{\rm us}(j)$ & $A_j$ &
$\rho_n^{(j)}$   & $\#\Omega_{\rm us}(j)$ & $\#\partial\Omega_{\rm us}(j)$ & $A_j$\\
\hline
$[0, 0, 0, 0, 0, 10]$ & $19$ & $18$ & $0.5$ & $[0, 0, 1, 0, 0, 9]$  & $166$ & $144$ & $12.5$ \\
\hline
$[0, 1, 0, 1, 0, 8]$  & $693$ &  $548$ & $21$ & $[0, 2, 0, 0, 1, 7]$ & $553$ &  $437$ & $25$ \\
\hline
$[1, 0, 0, 2, 0, 7]$  & $484$ &  $379$ & $25$ & $[0, 3, 0, 0, 0, 6]$  & $151$ &  $130$ & $9$\\
\hline
$[1, 1, 0, 1, 1, 6]$ & $2178$ &  $1574$ & $29$ & $[0, 0, 0, 3, 0, 6]$  & $118$ &  $100$ & $9$ \\
\hline
$[1, 2, 0, 1, 0, 5]$  & $1175$ &  $869$ & $29$ & $[0, 1, 0, 2, 1, 5]$ & $1035$ &  $756$ & $29$ \\
\hline
$[2, 0, 1, 0, 2, 5]$  & $1146$ &  $841$ & $29$ & $[0, 2, 0, 2, 0, 4]$  & $674$ &  $487$ & $17.5$\\
\hline
$[2, 1, 1, 0, 1, 4]$ & $2134$ & $1523$ & $29$ & $[3, 0, 0, 0, 3, 4]$  & $235$ &  $193$ & $9$ \\
\hline
$[1, 0, 1, 1, 2, 4]$  & $1899$ & $1333$ & $29$ & $[3, 1, 0, 0, 2, 3]$ & $800$ & $605$ & $25$ \\
\hline
$[1, 1, 1, 1, 1, 3]$  & $3909$ & $2590$ & $30$ & $[2, 0, 0, 1, 3, 3]$  & $727$ &  $541$ & $25$\\
\hline
$[3, 0, 2, 0, 0, 3]$  &  $497$ & $375$ & $25$ & $[0, 0, 2, 0, 3, 3]$  & $454$ & $339$ & $25$\\
\hline
$[2, 1, 0, 1, 2, 2]$  & $2108$ & $1436$ & $32$ & $[0, 1, 2, 0, 2, 2]$  & $1333$ & $918$ & $30$\\
\hline
$[4, 0, 1, 0, 1, 2]$  & $1116$ & $832$ & $21$ & $[1, 0, 1, 0, 4, 2]$  & $1008$ & $742$ & $21$\\
\hline
$[2, 0, 2, 1, 0, 2]$  & $1333$ & $918$ & $30$ & $[1, 1, 1, 0, 3, 1]$  & $2268$ & $1554$ & $32$\\
\hline
$[3, 0, 1, 1, 1, 1]$  & $2268$ & $1554$ & $32$ & $[1, 0, 3, 0, 1, 1]$  & $1308$ & $ 920$ & $30$\\
\hline
$[5, 0, 0, 1, 0, 1]$  & $413$ & $332$ & $13$ & $[0, 1, 0, 0, 5, 1]$  & $375$ & $299$ & $13$\\
\hline
$[0, 2, 0, 0, 4, 0]$  & $242 $ & $192$ & $17$ & $[2, 0, 2, 0, 2, 0]$  & $861$ & $598$ & $32$\\
\hline
$[0, 0, 0, 0, 6, 0]$  & $30$ & $27$ & $0.5$ & $[4, 0, 0, 2, 0, 0]$  &  $242$ & $192$ & $17$\\
\hline
$[0, 0, 4, 0, 0, 0]$  & $75$ & $66$ & $0.5$ & $[6, 0, 0, 0, 0, 0]$  & $37$ & $34$ & $0.5$
\end{tabular}
\label{table-EII}
\end{table}

We summarize the information in Table \ref{table-EII}. Moreover, $B=42$, which is attained at the $\frk$-type $[0, 0, 0, 0, 0, 2]$. The $E_{6(2)}$ case in Theorem \ref{thm-real-exceptional} follows.


\subsection{${\rm EIV}=E_{6(-26)}$}
\begin{figure}[H]
\centering
\scalebox{0.60}{\includegraphics{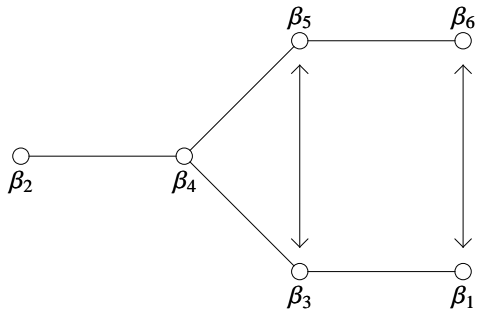}}
\caption{The Vogan diagram for EIV}
\label{Fig-EIV-Vogan}
\end{figure}

The Vogan diagram of this simple Lie algebra is presented in Fig.~\ref{Fig-EIV-Vogan}.
The positive root systems $(\Delta^+)^{(0)}(\frg, \frt_f)$ and $(\Delta^+)^{(0)}(\frk, \frt_f)$ are both of type $F_4$. The simple roots are $\alpha_1=\gamma_1:=\frac{1}{2}(\beta_1+\beta_6)$, $\alpha_2=\gamma_2:=\frac{1}{2}(\beta_3+\beta_5)$, $\alpha_3=\gamma_3:=\beta_4$, $\alpha_4=\gamma_4:=\beta_2$. The set $W(\frg, \frt_f)^1=\{e\}$. We have that $\rho_n^{(0)}=[1, 1, 0, 0]$ and that $\beta=[1, 0, 0, 0]$.

We compute that $\#\Omega_{\rm us}(0)=37$, $\#\partial\Omega_{\rm us}(0)=28$, $A_0=21$ and $B=37$, which is attained at the trivial $\frk$-type. The $E_{6(-26)}$ case in Theorem \ref{thm-real-exceptional} follows.

\subsection{${\rm EV}=E_{7(7)}$}
The Vogan diagram of this simple Lie algebra is obtained by painting the simple root $\alpha_2$ in Fig.~\ref{Fig-E7-Dynkin}. The positive roots system $\Delta^+(\frk, \frt_f)$ is of type $A_7$. Indeed, it has simple roots
$$
\gamma_1:=\alpha_1; \quad \gamma_i:=\alpha_{i+1}, \, 2\leq i\leq 6; \quad \gamma_7:=\alpha_1+2\alpha_2+2\alpha_3+3\alpha_4+2\alpha_5+\alpha_6.
$$
Note that $|W(\frg, \frt_f)^1|=72$ and $\beta=[0,0,0, 1, 0, 0, 0]$.

We calculate that $\max\{A_j\mid 0\leq j\leq 71\}=78$.
Moreover, $B=\frac{231}{2}$, which is attained at the trivial $\frk$-type. The $E_{7(7)}$ case in Theorem \ref{thm-real-exceptional} follows.

\subsection{${\rm EVI}=E_{7(-5)}$}
The Vogan diagram of this simple Lie algebra is obtained by painting the simple root $\alpha_1$ in Fig.~\ref{Fig-E7-Dynkin}.  The positive roots system $\Delta^+(\frk, \frt_f)$ is of type $D_6\times A_1$. Indeed, it has simple roots
$$
\gamma_i:=\alpha_{8-i}, \, 1\leq i\leq 4; \, \gamma_5=\alpha_2; \, \gamma_6=\alpha_3;
\, \gamma_7:=2\alpha_1+2\alpha_2+3\alpha_3+4\alpha_4+3\alpha_5+2\alpha_6+\alpha_7.
$$
Note that $|W(\frg, \frt_f)^1|=63$ and $\beta=[0, 0, 0, 0, 0, 1, 1]$.

We compute that $\max\{A_j\mid 0\leq j\leq 62\}=82$. Moreover, $B=\frac{227}{2}$, which is attained at the $\frk$-type $\mu=[0, 0, 0, 0, 0, 0, 4]$. Indeed, ${\rm MP}(\mu)=\|\mu+\beta\|_{\rm spin}=\frac{231}{2}$, while $\|\mu\|_{\rm lambda}=2$.  The $E_{7(-5)}$ case in Theorem \ref{thm-real-exceptional} follows.


\subsection{${\rm EVIII}=E_{8(8)}$}
The Vogan diagram of this simple Lie algebra is obtained by painting the simple root $\alpha_1$ in Fig.~\ref{Fig-E8-Dynkin}.  The positive roots system $\Delta^+(\frk, \frt_f)$ is of type $D_8$. Indeed, it has simple roots
$$
\gamma_1:=2\alpha_1+2\alpha_2+3\alpha_3+4\alpha_4+3\alpha_5+2\alpha_6+\alpha_7;
\,\gamma_i:=\alpha_{10-i}, \, 2\leq i\leq 6; \, \gamma_7=\alpha_2; \, \gamma_8=\alpha_3.
$$
Note that $|W(\frg, \frt_f)^1|=135$ and $\beta=[0, 0, 0, 0, 0, 0, 1, 0]$.

We compute that $\max\{A_j\mid 0\leq j\leq 134\}=\frac{411}{2}$. Moreover, $B=380$, which is attained at the trivial $\frk$-type. The $E_{8(8)}$ case in Theorem \ref{thm-real-exceptional} follows.


\subsection{${\rm EIX}=E_{8(-24)}$}

The Vogan diagram of this simple Lie algebra is obtained by painting the simple root $\alpha_8$ in Fig.~\ref{Fig-E8-Dynkin}.  The positive roots system $\Delta^+(\frk, \frt_f)$ is of type $E_7\times A_1$. Indeed, it has simple roots
$$
\gamma_i:=\alpha_{i}, \, 1\leq i\leq 7; \quad \gamma_8:=2\alpha_1+3\alpha_2+4\alpha_3+6\alpha_4+5\alpha_5+4\alpha_6+3\alpha_7+2\alpha_8.
$$
Note that $|W(\frg, \frt_f)^1|=120$ and that $\beta=[0, 0, 0, 0, 0, 0, 1, 1]$.

We compute that $\max\{A_j\mid 0\leq j\leq 119\}=212$.
Moreover, $B=\frac{723}{2}$, which is attained at the $\frk$-type $[0, 0, 0, 0, 0, 0, 0, 8]$. The $E_{8(-24)}$ case in Theorem \ref{thm-real-exceptional} follows.

\section{Hermitian symmetric case}\label{sec-HS}

In this section, we consider the case that $(G, K)$ is a Hermitian symmetric pair. Then the $K$-types of an infinite-dimensional $(\frg, K)$ module $\pi$ may \emph{not} be a union of Vogan pencils.
Therefore, if $\pi$ is a $(\frg, K)$ module with a lowest $K$-type $E_{\mu}$ and with infinitesimal character \eqref{inf-char}, we can only conclude that
\begin{equation}\label{nu-bound-HS}
\|\nu\|^2\leq \|\mu\|_{\rm spin}^2 - \|\mu\|_{\rm lambda}^2.
\end{equation}
Thus our job is to get the maximum of the RHS of \eqref{nu-bound-HS} over all the $K$-types. Our main result of this section is summarized as follows.

\begin{thm}\label{thm-real-exceptional-HS}
Let $G$ be the Lie group \emph{\texttt{E6\_h}} or \emph{\texttt{E7\_h}} in \emph{\texttt{atlas}} \cite{ALTV, At}. Let $\pi$ be any irreducible unitary $(\frg, K)$ module whose infinitesimal character is given by \eqref{inf-char}. Then we have the following bounds for $\|\nu\|$:
\begin{itemize}
\item[$\bullet$] \emph{\texttt{E6\_h}:} $\|\nu\|\leq \sqrt{73}$, while $\|\rho(G)\|=\sqrt{78}$.
\item[$\bullet$] \emph{\texttt{E7\_h}:} $\|\nu\|\leq \sqrt{\frac{371}{2}}$, while $\|\rho(G)\|=\sqrt{\frac{399}{2}}$.
\end{itemize}
In each case, the bound is attained at the trivial representation.
\end{thm}

Let $\frt_f$ be the fundamental Cartan subalgebra of $\frg$. We may and we will identify $\frt_f$ with $\frt_f^*$ via the Killing form $B(\cdot, \cdot)$. Let $\bbR\zeta$ be the one-dimensional center of $\frk_f$, and let $\frt_f^-$ be the orthogonal complement of $\bbR\zeta$ in $\frt_f$ under $B(\cdot, \cdot)$.

Let $|W(\frg, \frt_f)^1|=s$, and enumerate its elements as $w^{(j)}$ for $0\leq j\leq s-1$. We adopt the notation $\Omega(j)$, $\Omega_{\rm us}(j)$, $\partial\Omega_{\rm us}(j)$, $A_j$ etc as in Section \ref{sec-framework}, with the discrepancy that instead of $\frk$-types, now we use $K$-types. Put
\begin{equation}\label{bound-B-HS}
B=\max_{\mu}\left\{ \|\mu\|^2_{\rm spin} - \|\mu\|^2_{\rm lambda}  \right\},
\end{equation}
where $\mu$ runs over all the $K$-types in $\Omega_{\rm us}(j)\setminus \partial\Omega_{\rm us}(j)$. Then $\|\nu\|^2$ is upper bounded by $\max\{A_0, A_1, \dots, A_{s-1}, B\}$.

\subsection{\texttt{E6\_h}}
This non-compact simple linear Lie group in \texttt{atlas} has center $\bbZ/3\bbZ$ and Lie algebra ${\rm EIII}=E_{6(-14)}$. Its Vogan diagram is obtained from Fig.~\ref{Fig-E6-Dynkin} by painting the simple root $\alpha_6$. Let $\zeta_1, \dots, \zeta_6:=\zeta$ be the fundamental weights corresponding to the simple roots $\alpha_1, \dots, \alpha_6$ of $\Delta^+(\frg, \frt_f)$. Put $\gamma_i=\alpha_i$ for $1\leq i\leq 5$. Then $\gamma_1, \dots, \gamma_5$ are the simple roots of $\Delta^+(\frk, \frt_f^-)$, which is of type $A_5$. Let $\varpi_1, \dots, \varpi_5\in (\frt_f^-)^*$ be the corresponding fundamental weights. For $a, b, c, d, e\in\bbZ_{\geq 0}$ and  $f\in\bbZ$, let $[a, b, c, d, e, f]$ stand for the vector $a \varpi_1+b \varpi_2 + c \varpi_3+d \varpi_4 +e \varpi_5+ \frac{f}{4}\zeta$. Then ${[a, b, c, d, e, f]}$ is the highest weight of a $K$-type if and only if
\begin{equation}\label{EIII-K-type}
-\frac{3a}{4}  -\frac{5b}{4} -\frac{3c}{2} -d -\frac{e}{2} +\frac{f}{4} \in\bbZ.
\end{equation}

We have that $|W(\frg, \frt_f)^1|=27$. We compute that $\max\{A_j\mid 0\leq j\leq 26\}=33$.
Moreover, $B=73$, which is attained at the trivial $K$-type. The \texttt{E6\_h} case in Theorem \ref{thm-real-exceptional-HS} follows.

\begin{example}\label{exam-EIII}
Take the index $j$ so that $\rho_n^{(j)}=[1, 0, 1, 0, 1, 3]$ and consider the $K$-type $\mu=[a, b, c, d, e, f]$.
Then $\mu\in\Omega(j)$ if and only if
\begin{align*}
b+2e\leq a+2c+f, a+f\leq b+2c+2e+8, b+2c\leq a+2e+f, 2e+f\leq 3a+b+2c+8.
\end{align*}
On the other hand,
\begin{align*}
&w^{(j)}\zeta_1=[0, 1, 0, 0, 0, 1], \quad w^{(j)}\zeta_2=[0, 0, 0, 1, 0, 0], \quad w^{(j)}\zeta_3=[0, 0, 1, 0, 0, 2], \\
&w^{(j)}\zeta_4=[0, 0, 1, 0, 1, 0], \quad w^{(j)}\zeta_5=[1, 0, 0, 0, 1, 1], \quad w^{(j)}\zeta_6=[1, 0, 0, 0, 0, -1].
\end{align*}
We compute that the non-decreasable $K$-types in $\Omega(j)$ are
\begin{align*}
[0, 0, 0, 0, 0, 0], [0, 0, 0, 0, 0, 4], [0, 0, 0, 0, 0, 8], [0, 0, 0,
0, 1, 2], [0, 0, 0, 0, 1, 6], [0, 0, 0, 0, 2, 4].
\end{align*}
They are all u-small. \hfill\qed
\end{example}

\subsection{\texttt{E7\_h}}

This non-compact simple linear Lie group in \texttt{atlas} has center $\bbZ/2\bbZ$ and Lie algebra ${\rm EVII}=E_{7(-25)}$. Its Vogan diagram is obtained from Fig.~\ref{Fig-E7-Dynkin} by painting the simple root $\alpha_7$. Let $\zeta_1, \dots, \zeta_7:=\zeta$ be the fundamental weights corresponding to the simple roots $\alpha_1, \dots, \alpha_7$ of $\Delta^+(\frg, \frt_f)$. Put $\gamma_i=\alpha_i$ for $1\leq i\leq 6$. Then $\gamma_1, \dots, \gamma_6$ are the simple roots of $\Delta^+(\frk, \frt_f^-)$, which is of type $E_6$. Let $\varpi_1, \dots, \varpi_6\in (\frt_f^-)^*$ be the corresponding fundamental weights. For $a, b, c, d, e, f\in\bbZ_{\geq 0}$ and  $g\in\bbZ$, let $[a, b, c, d, e, f, g]$ stand for the vector $a \varpi_1+b \varpi_2 + c \varpi_3+d \varpi_4 +e \varpi_5+ f\varpi_6 + \frac{g}{3}\zeta$. Then $[a, b, c, d, e, f, g]$ is the highest weight of a $K$-type if and only if
\begin{equation}\label{EIII-K-type}
-\frac{2 a}{3}-b-\frac{4 c}{3}-2 d-\frac{5 e}{3}-\frac{4 f}{3}+ \frac{g}{3} \in\bbZ.
\end{equation}

We have that $|W(\frg, \frt_f)^1|=56$. We compute that $\max\{A_j\mid 0\leq j\leq 55\}=79$. Moreover, $B=\frac{371}{2}$, which is attained at the trivial $K$-type. The \texttt{E7\_h} case in Theorem \ref{thm-real-exceptional-HS} follows.

\medskip
\centerline{\scshape Funding}
Dong is supported by the National Natural Science Foundation of China (grant 12171344, 2022-2025).

\medskip
\centerline{\scshape Acknowledgements}
I thank Daniel Wong for helpful discussions.
I thank an anonymous referee for giving us expertise suggestions.

\end{document}